\documentclass[11pt]{amsart}
\usepackage{amsmath, amsthm, latexsym, amssymb, amsfonts}
\usepackage[pdftex]{graphicx,color}

\newenvironment{pf}{\proof[\proofname]}{\endproof}
\theoremstyle{plain}
\newtheorem{Th}{Theorem}[section]
\newtheorem{Cor}[Th]{Corollary}
\newtheorem{Prop}[Th]{Proposition}
\newtheorem{Lem}[Th]{Lemma}
\numberwithin{equation}{section} \theoremstyle{definition}
\newtheorem{Rem}[Th]{Remark}

\newtheorem{Ex}[Th]{Example}

\newtheorem{Prob}{Problem} 

\newtheorem{Def}[Th]{Definition}

\usepackage{mathtools}

\newcommand{\cal}[1]{\mathcal{#1}}
\newcommand{\La}{\mathbb L}

\newcommand{\B}{\mathbb B}
\newcommand{\Q}{\mathbb Q}
\newcommand{\Z}{\mathbb Z}
\newcommand{\R}{\mathbb R}
\newcommand{\N}{\mathbb N}

\newcommand{\la}{\lambda}
\newcommand{\al}{\alpha}
\newcommand{\D}{\Delta}

\newcommand{\cL}{\cal L}

\newcommand{\cZ}{\cal Z}

\newcommand{\Vol}{\operatorname{Vol}}
\newcommand{\chull}{\operatorname{Conv}}
\newcommand{\ve}{\varepsilon}

\newcommand{\rs}[1]{Section~\ref{S:#1}}
\newcommand{\rl}[1]{Lemma~\ref{L:#1}}
\newcommand{\rp}[1]{Proposition~\ref{P:#1}}

\newcommand{\rr}[1]{Remark~\ref{R:#1}}

\newcommand{\re}[1]{(\ref{e:#1})}
\newcommand{\rc}[1]{Corollary~\ref{C:#1}}
\newcommand{\rt}[1] {Theorem~\ref{T:#1}}
\newcommand{\rd}[1]{Definition~\ref{D:#1}}
\newcommand{\rf}[1]{Figure~\ref{F:#1}}

\begin{document}
\author{Ivan Soprunov} 
\address{Mathematics Department, Cleveland State University, Cleveland, OH 44115}
\email{i.soprunov@csuohio.edu}
\author{Jenya Soprunova}
\address{Department of Mathematical Sciences,
Kent State University, Kent, OH 44242}
\email{soprunova@math.kent.edu}

\thanks{The first author is partially supported by NSA Grant H98230-13-1-0279}

\title{Eventual Quasi-linearity of the Minkowski length}

\begin{abstract}  The Minkowski length of a lattice polytope $P$
is a natural generalization of the lattice diameter of $P$.
It can be defined as the largest number of lattice segments whose Minkowski sum is contained in $P$.
The famous Ehrhart theorem states that the number of lattice points in the positive integer dilates $tP$
of a lattice polytope $P$ behaves polynomially in $t\in\N$.
In this paper we prove that for any lattice polytope $P$,
the Minkowski length of $tP$ for $t\in\N$ is eventually a quasi-polynomial with linear constituents.
We also give a formula for the Minkowski length of coordinates boxes, degree one polytopes, 
and dilates of unimodular simplices. In addition, we give a new bound for the Minkowski length of lattice polygons
and show that the Minkowski length of a lattice triangle coincides with its lattice diameter.
\end{abstract}

\maketitle

\section*{Introduction}
Let $P$ be a $d$-dimensional lattice polytope in $\R^d$. Recall that the lattice
diameter $\ell(P)$ is defined as one less than the largest number of collinear lattice points in $P$.
The Minkowski length is a natural extension of this notion. For any
$1\leq n\leq d$, let  $L_n(P)$ be the largest number of lattice polytopes of positive dimension whose Minkowski sum 
is at most $n$-dimensional and is contained in $P$. We call $L_n(P)$ the $n$-th Minkowski length of $P$,
and $L(P)=L_d(P)$ simply the Minkowski length of $P$.  Note $L_1(P)$ coincides with the lattice diameter $\ell(P)$,
as in this case the Minkowski summands are collinear lattice segments. It is not hard to show (see the discussion after \rd{maximal})
that $L_n(P)$ is the largest number of lattice segments whose Minkowski sum is at most $n$-dimensional 
and is contained in $P$.

The Minkowski length $L(P)$ of a lattice polytope $P\subset\R^d$ was first introduced in \cite{SS1}
in relation to studying parameters of toric surface codes. Every lattice polytope $P$ defines
a space $\cL(P)$ of Laurent polynomials (over some field) whose monomials have exponent vectors lying in $P$.
Such spaces naturally appear in the theory of toric varieties. 
The algebraic interpretation of the Minkowski length is the following: $L(P)$ is the largest number
of irreducible factors a polynomial $f\in\cL(P)$ may have. This information is particularly important 
when one studies zeroes of polynomials in $\cL(P)$, see \cite{LSc, Lit11, SS1, SS2, UmVe}.  
 A number of results concerning $L(P)$ appeared in \cite{BGSW, SS1, Whit}.

Let $tP=\{tx\in\R^d\ | x\in P\}$ be the dilate of $P$ by a positive integer factor $t$.
 The main result of this paper explains 
 the behavior of $L_n(tP)$ as a function of the scaling factor $t\in\N$ in the spirit of the Ehrhart theory.
 In \rt{quasi-lin} we prove that for any lattice polytope $P$ the function $L_n(tP)$ is eventually quasi-polynomial
 with linear constituents (we say ``quasi-linear" for short), which
 contributes positively to the ``ubiquitousness of quasi-polynomials" phenomenon declared by Kevin Woods \cite{Kevin}.
 For an introduction to the Ehrhart theory we refer the reader to the wonderful book 
 by M.~Beck and S.~Robins \cite{BR}.
 
To prove eventual quasi-linearity of the Minkowski lengths $L_n(P)$
we define and study their rational counterparts: 
a sequence of rational numbers $\la_1(P)\leq\dots\leq \la_d(P)=\la(P)$ associated with $P$. Here $\la_1(P)$ is the rational diameter of 
$P$ and $\la_n(P)$ is the ``asymptotic" Minkowski length, i.e. $\la_n(P)=\lim_{t\to\infty}L_n(tP)/t$.
In \rt{period} we prove that $\la_n(P)=L_n(kP)/k$ for some $k\in\N$.

Although an algorithm for computing $L(P)$ was presented in dimensions two and three (see \cite{BGSW, SS1}),
 there have been no explicit formulas for $L(P)$ even for simplices. Here we prove that
 $L(t\D)=t$ for any unimodular simplex $\D$ and any $t\in\N$.
This result allowed us to write explicit answers for $L(P)$ for other classes of polytopes
such as coordinate boxes and polytopes of degree one (see \rc{simplex} and examples afterwards). 
In \rs{dim2} we prove that for lattice triangles the Minkowski length coincides with the lattice diameter.  
The final part of the paper contains some examples and open questions.

\section{Preliminaries}

We start with some standard terminology from geometric combinatorics. 
A polytope $P\subset\R^d$ is called {\it lattice} (resp. {\it rational}\,) if its vertices have integer (resp. rational) coordinates. A vector $v\in\Z^d$ is called {\it primitive} if the greatest common divisor of its coordinates is 1.
A lattice segment is called {\it primitive} if it contains exactly two lattice points. 
A $d$-dimensional simplex is called {\it unimodular} if its vertices affinely generate the lattice $\Z^d$.

Given a lattice polytope $P\subset\R^d$, denote by $\Vol_d(P)$ the Euclidean $d$-dimensional 
volume of $P$. Note that the $d$-dimensional volume of any parallelepiped formed by a basis of $\Z^d$
equals 1. More generally, suppose $P$ is contained in an $n$-dimensional rational affine subspace $a+H$
for a rational linear subspace $H\subset\R^d$ and $a\in\Q^d$.
Denote by $\Vol_{n}(P)$ the $n$-dimensional volume of $P$ normalized such that the 
$n$-dimensional volume of any parallelepiped formed by a basis of the lattice $H\cap\Z^d$
equals 1.

\subsection{Minkowski length}

Let $P$ and $Q$ be convex polytopes in $\R^d$. Their {\it Minkowski sum} is the set
$$P+Q=\{p+q\in\R^d\ |\ p\in P,\ q\in Q\},$$
which is again a convex polytope.

\begin{Def}\label{D:maximal} Let $P$ be a lattice polytope in $\R^d$.
Define the {\it Minkowski length} $L=L(P)$ of $P$ to be the largest number
of lattice polytopes $Q_1,\dots, Q_L$ of positive dimension  whose 
Minkowski sum is contained in $P$. 
Any such sum $Q_1+\dots+Q_L$ is called a {\it maximal decomposition in P}. 
\end{Def}

We refer the reader to \cite{BGSW} for examples illustrating this definition.
It is clear from the definition that $L(P)$ is monotone with respect to inclusion: $L(P)\leq L(Q)$ if $P\subseteq Q$,
and is superadditive with respect to the Minkowski sum:  $L(P+Q)\geq L(P)+L(Q)$.
Also, $L(P)$ is invariant under unimodular transformations (isomorphisms of the lattice $\Z^d$).

There is a natural partial order on the set of maximal decompositions in $P$, as defined in \cite{BGSW}.
Namely, we say that 
$$Q_1+\cdots+Q_L\leq P_1+\cdots+P_L$$ 
if  $Q_1+\cdots+Q_L$ is contained in $P_1+\cdots+P_L$ after a possible lattice translation. 
Minimal elements with respect to this partial order are called {\it smallest maximal decompositions}. 
Clearly, every smallest maximal decomposition is  the Minkowski sum of $L$ lattice segments, i.e.
is a lattice zonotope.

Any lattice (resp. rational) zonotope $Z$ can be written in the form
\begin{equation}\label{e:zonotope}
Z=a+\al_1[0,v_1]+\dots+\al_m[0,v_m],
\end{equation}
for some $m\in\N$, distinct primitive vectors $v_i\in\Z^d$, positive integer (resp. rational) numbers $\al_i$, and a
lattice (resp. rational) point $a\in \R^d$. In this case, we set 
$$|Z|:=\al_1+\dots+\al_m.$$

The following result from  \cite{BGSW} gives a universal bound for the number of distinct summands in
a smallest maximal decomposition.

\begin{Prop}\label{P:bound} \cite{BGSW}
Let $P\subset\R^d$ be a lattice polytope. Then every smallest maximal decomposition in $P$
has at most $2^d-1$ distinct summands.
\end{Prop}

\subsection{Quasi-polynomials} Here we recall the definition of a quasi-polynomial function.

\begin{Def}
A function $f:\N\to \Q$ is called a {\it quasi-polynomial} if there exist $k\in\N$
and polynomials $p_0,\dots, p_{k-1}\in\Q[t]$, called the {\it constituents} of $f$, such that 
$$f(t)=p_r(t)\quad\text{whenever }t\equiv r\!\!\!\mod k, \text{ for }0\leq r\leq k.$$
The smallest such $k$ is called the {\it period} of $f$. If all the constituents of $f$ are linear we say that $f$ is {\it quasi-linear}. Finally,
we say  $f:\N\to \Q$ is {\it eventually quasi-linear} if $f(t)$ coincides with a
quasi-linear function for all large enough $t$.
\end{Def}

\begin{Ex}
The function$f(t)= 3\left\lfloor\frac{t}{3}\right\rfloor+4$ is quasi-linear with period $k=3$, where
$\lfloor x\rfloor$ denotes the floor of $x$. Indeed, 
$$f(t)=\begin{cases}t+4,& \text{if }\  t\equiv 0\!\!\mod 3\\
                                    t+3,& \text{if }\ t\equiv 1\!\!\mod 3\\
                                    t+2,& \text{if }\ t\equiv 2\!\!\mod 3.
                                    \end{cases}$$
\end{Ex}

\section{Main Theorems}

Before we prove our main result about eventual quasi-linearity of the Minkowski length, we
will look at some instances when it is, in fact, linear. The simplest such example is
when $P=\D$, a unimodular $d$-simplex.

\begin{Th}\label{T:box} 
Let $\D$ be a unimodular $d$-simplex and $t\in\N$. Then $$L(t\D)=tL(\D)=t.$$
\end{Th}

\begin{pf} After a unimodular transformation we may assume
that $\D$ is the standard $d$-simplex, i.e. the convex hull of $\{0,e_1,\dots,e_d\}$
where $e_i$ are the standard basis vectors. First, it is easy to see that $L(\D)=1$. 
Also, $L(t\D)\geq  t$ as $t\D$ contains
the Minkowski sum of $t$ lattice segments $[0,e_1]$.

We prove the converse by induction on $d$. The case $d=1$ is trivial.
Denote $L=L(t\D)$. Let $Z$ be a smallest maximal decomposition in $t\D$ and let
$a\in Z$ be a vertex with the smallest sum of the coordinates, which we denote by $\al$. 
We have
$$Z=a+[0,v_1]+\dots+[0,v_L],$$
where $v_i\in\Z^d$ are primitive, not necessarily distinct vectors.
Note that the sum of the coordinates of each $v_i$ is non-negative, by the choice of $a$.
Suppose the first $k$ of the vectors $v_1,\dots, v_L$  have the sum of the coordinates equal zero, for $0\leq k\leq L$.
Then the subzonotope 
$$Z'=a+[0,v_1]+\dots+[0,v_k]$$
is contained in $\al\D'$, where $\D'$ is the facet of $\D$ whose points have the sum of the coordinates equal to 1.
This implies $k\leq L(\al\D')$.  By induction $L(\al\D')=\al$, hence, $k\leq\al$. 
Now the point $v=a+v_{k+1}+\dots+v_L$ lies in $Z$ and
has the sum of the coordinates at least $\al+L-k\geq L$. On the other hand, $v$ lies in $t\D$,
so its sum of the coordinates is at most~$t$. Therefore, $L\leq t$. 
\end{pf}

\begin{Cor}\label{C:simplex} 
Let $P$ be a lattice polytope contained in $\al\D$ for some unimodular simplex $\D$ and $\al\in\N$. If $P$ contains 
the Minkowski sum of $\al$ lattice segments then $L(P)=\al$. Consequently, $L(tP)=tL(P)$.
\end{Cor}

\begin{Ex}
Let $\Pi$ be a lattice coordinate box in $\R^d$, i.e. 
$\Pi=[0,\al_1e_1]\times\dots\times [0,\al_de_d]$, where $e_i$ are the standard basis vectors and $\al_i$ are
non-negative integers. Clearly, $\Pi=\al_1[0,e_1]+\dots+\al_d[0,e_d]$ and $\Pi$ is contained in 
$(\al_1+\dots+\al_d)\D_d$, where $\D_d$ is the standard $d$-simplex. Therefore,
$$L(\Pi)=\al_1+\dots+\al_d.$$
We also have $L(t\Pi)=tL(\Pi)$.
\end{Ex}

\begin{Ex} According to a result of Batyrev and Nill \cite{BN},  
a $d$-dimensional polytope $P$ has {\it degree one} if and only if $P$ is either
\begin{enumerate}
\item the $d-2$ iterated pyramid over the triangle $\D_2=\chull\{(0,0),(2,0),(0,2)\}$, or
\item the $d-n$ iterated pyramid over a Lawrence prism $Q$ 
defined by a sequence of integers $0<h_1\leq\dots\leq h_{n}$:
\begin{equation}\label{e:Lawrence}
Q=\chull\{0,e_1,\dots, e_{n-1}, e_1+h_1e_n,\dots, e_{n-1}+h_{n-1}e_n,h_ne_n\}\subset\R^n.
\end{equation}
\end{enumerate}
 \rc{simplex} implies that in the first case $L(P)=2$ since $P$ contains the segment $[0,2e_1]$ and $P\subset 2\D_d$.
In the second case
$$L(P)=\begin{cases}h_n, &\text{ if }h_{n-1}<h_n\\ h_n+1, &\text{ if }h_{n-1}=h_n.\end{cases}$$ 
Indeed, if $h_{n-1}<h_n$ then $P\subset h_n\D_d$ and $P$ contains the segment $[0,h_ne_n]$.
If $h_{n-1}=h_n$ then $P\subset (h_n+1)\D_d$ and $P$ contains the rectangle $[0,e_{n-1}]+[0,h_ne_n]$.
\end{Ex}

\subsection{Rational Minkowski length}\label{S:lambda}
Let $P$ be an arbitrary lattice polytope in $\R^d$. 
The following is a generalization of \rd{maximal}.

\begin{Def}\label{D:genmaximal} Let $P$ be a lattice polytope in $\R^d$.
Define the {\it $n$-th Minkowski length} $L=L_n(P)$ of $P$ to be the largest number
of lattice polytopes $Q_1,\dots, Q_L$ of positive dimension  whose 
Minkowski sum is at most $n$-dimensional and is contained in $P$.
\end{Def}

Clearly, $L_1(P)\leq \dots \leq L_{d-1}(P)\leq L_d(P)=L(P)$. Note $L_1(P)$
coincides with the  {\it lattice diameter}  $\ell(P)$, which is defined as one less than
the largest number of collinear lattice points in $P$.

\begin{Ex} Let $\square$ be the unit square in $\R^2$. Then 
$$L_1(\square)=1\quad\text{ and }\quad L_2(\square)=L(\square)=2.$$
For any unimodular $d$-simplex $\D$ and any $t\in\N$ we have $L_1(t\D)\geq t$
as $t\D$ contains the segment $[0,te_1]$. By \rt{box}, $L(t\D)=t$, hence
$$L_1(t\D)=\dots=L_d(t\D)=L(t\D)=t.$$
If $P$ is the $d-n$ iterated pyramid over a Lawrence prism $Q$ as in \re{Lawrence} then
$L_1(P)=h_n$ and $L_2(P)=\dots=L_d(P)=L(P)$.
\end{Ex}

The following is a rational analog of the Minkowski length.

\begin{Def}\label{D:lambda} The number 
$$\la(P)=\sup_{t\in\N}\frac{L(tP)}{t}$$
is called the {\it rational Minkowski length} of $P$. More generally,
the {\it $n$-th rational  Minkowski length of $P$} is
$$\la_n(P)=\sup_{t\in\N}\frac{L_n(tP)}{t}$$
\end{Def}

The following proposition asserts that the numbers $\la_n(P)$ are well-defined.

\begin{Prop} For any $t\in\N$ we have $L_n(tP)\leq t\al$, where $\al\in\N$ is such that
$P\subseteq \al\D$ for a unimodular simplex $\D$.
\end{Prop}

\begin{pf} It is enough to consider the case $n=d$. Then it is immediate from \rt{box}:
$L(tP)\leq L(t\al\D)=t\al$.
\end{pf}

It follows from the definition that $\la_1(P)\leq\dots\leq\la_{d-1}(P)\leq\la_d(P)=\la(P)$.

\begin{Rem} As $L_n(tP)$ satisfies the superadditivity property, the supremum in the above 
definition may be replaced with the limit, by Fekete's lemma. We will not be using
this result in our further discussion.
\end{Rem}

\begin{Def}\label{D:diameter} 
Let $K\subset\R^d$ be a rational polytope. For any primitive vector $v\in\Z^d$ define
$s_v(K)$ to be the largest rational number $s$ such that the segment $[0,sv]$ is contained in $K$
after a translation by a rational vector. The {\it rational diameter}  $s(K)$ is 
the maximum of $s_v(K)$ over all primitive $v\in\Z^d$.
\end{Def}

It is not hard to see that $s(P)=\la_1(P)$ for any lattice polytope $P$. Indeed,
for any $t\in\N$, the polytope $P$ contains a segment $a+[0,sv]$ for some $a\in\Q^d$, primitive $v\in\Z^d$,
and $s=L_1(tP)/t$. Thus, $\la_1(P)\leq s(P)$. Conversely, if $a+[0,s(P)v]$ is contained in $P$ for
some $a\in\Q^d$ and primitive $v\in\Z^d$ then there exists $t\in\N$ such that $ta+[0,ts(P)v]$
is a lattice segment contained in $tP$, i.e. $s(P)\leq L_1(tP)/t\leq \la_1(P)$.
As a corollary, we obtain $L_1(P)=\lfloor \la_1(P)\rfloor$, as $\ell(P)=\lfloor s(P)\rfloor$.

In our main theorem below (\rt{period}) we show that $\la(P)$ as well as all $\la_n(P)$ are, in fact,  rational numbers.
First, we need a few lemmas.

\begin{Lem}\label{L:uepsilon}
Let $K$ be a convex body in $\R^d$ and fix $\ve>0$. Then the set
$$U_\ve(K)=\{v\in\Z^d\ |\ v \text{ primitive},\, s_v(K)\geq\ve \}$$
is finite.
\end{Lem}

\begin{pf} First, note that if 
$K\subseteq K'$ then $s_v(K)\leq s_v(K')$, 
and $s_v(\al K)=\al s_v(K)$ for $\al\in\Q$. Thus it is enough
to prove the statement for $K=\B$, the $d$-dimensional unit ball. 
Let $v\in\Z^d$ be primitive. By definition $s_v(\B)$ is the number $s\in\Q$ such that 
$\|sv\|=2$, where $\|\ \|$ is the usual Euclidean norm. It follows that $s_v(\B)\geq\ve$
if and only if $\|v\|\leq 2/\ve$. In other words, $U_\ve(\B)$ is a lattice set contained in 
the ball of radius $2/\ve$ and so is finite.
\end{pf}

\begin{Lem}\label{L:volume}
Let $Z=a+\al_1[0,v_1]+\dots+\al_m[0,v_m]$ be a smallest maximal decomposition in $P$
and $n=\dim Z$. Then for any $1\leq i_1<\dots<i_n\leq m$ the $n$-dimensional volume 
of the parallelepiped formed by $v_{i_1},\dots, v_{i_n}$ is no greater than~$n^d$.
\end{Lem}

\begin{pf} We may assume that $v_{i_1},\dots, v_{i_n}$ are linearly independent.
Let $a+H$ be the affine span of $Z$ and let $\La=H\cap\Z^d$
be the corresponding lattice of rank $n$.
It is well-known that the $n$-dimensional volume of the parallelepiped formed by $n$ linearly independent lattice
vectors $w_1,\dots, w_n\in\La$ equals the number of lattice points in the half-open parallelepiped
$$\{\la_1w_1+\dots+\la_nw_n\ |\ 0\leq\la_i<1 \text{ for }1\leq i\leq n\},$$
which is less than the number of lattice points in the closed parallelepiped.

Let $\{v_{i_1},\dots, v_{i_n}\}$ be a subset of the set of vectors appearing in the decomposition $Z$.
We claim that the parallelepiped $\Pi$ they form has at most $n^d$ lattice points. Indeed,
consider the image of $\Pi\cap\Z^d$ in $(\Z/n\Z)^d$. If there are two lattice points in $\Pi$
congruent mod $(n\Z)^d$ then the lattice segment $E$ containing them lies in $\Pi$ and has lattice length $n$.
Therefore, if we replace $[0,v_{i_1}]+\dots+[0,v_{i_n}]$ in the decomposition $Z$ with $E$
we obtain a maximal decomposition $Z'$ properly contained in $Z$. This contradicts the fact
that $Z$ is a smallest maximal decomposition. This shows that all lattice points of $\Pi$
are different in $(\Z/n\Z)^d$, i.e. their number cannot exceed~$n^d$.
\end{pf}

\begin{Lem}\label{L:finite} Let $B=\{u_1,\dots,u_n\}$ be a basis for a rational linear subspace $H\subseteq\R^d$,
and fix a constant $N$. Let $\Vol_n(u_1,\dots, \underset{i}{v},\dots,u_n)$ denote the
$n$-dimensional volume  of the parallelepiped formed by  $u_1,\dots,u_n$ where
$u_i$ is replaced with a vector $v$. 
Then the set
$$V(B)=\{v\in H\cap\Z^d\ |\ \Vol_n(u_1,\dots, \underset{i}{v},\dots,u_n)\leq N, \text{ for all }1\leq i\leq n\}$$
is finite.
\end{Lem}

\begin{pf} 
Fix a coordinate system in $H$ by choosing a basis for the lattice $H\cap\Z^d$.
Write $v=x_1u_1+\dots+x_nu_n$ for some $x_i\in\R$. Then by Cramer's rule
$$|x_i|=\frac{\Vol_n(u_1,\dots, \underset{i}{v},\dots,u_n)}{\Vol_n(u_1,\dots,u_n)}\leq \frac{N}{\Vol_n(u_1,\dots,u_n)}=:c_i.$$
Therefore, the set $V(B)$ is contained in the set of lattice points of the parallelepiped formed by $\{\pm c_1u_1,\dots,\pm c_nu_n\}$,
and, hence, is finite.
\end{pf}

\begin{Lem}\label{L:allzonotopes} Let $P\subset\R^d$ be a lattice polytope.
Fix an ordered collection of primitive vectors ${\bf v}=(v_1,\dots, v_m)\in(\Z^d)^m$. Then the set of zonotopes
\begin{equation}
\cZ({\bf v})=\{Z=a+\al_1[0,v_1]+\dots+\al_m[0,v_m]\ |\ a\in\R^d, \al_i\in\R_{\geq 0}, Z\subseteq P\}\nonumber
\end{equation}
is a rational polytope in $\R^{d+m}$. The function $|\cdot|:\cZ({\bf v})\to\R$, $Z\mapsto |Z|$ is an integer linear 
function on $\cZ({\bf v})$.
\end{Lem}

\begin{pf} With every such zonotope $Z$ we associate a point $z=(a,\al_1,\dots,\al_m)\in\R^{d+m}$.
Note that $Z$ is the convex hull of the following set of  $2^m$ points in $\R^d$:
$$K=\Big\{a+\sum_{i\in I}\al_iv_i\ |\ I\subseteq \{1,\dots, m\}\Big\}.$$
Clearly $Z\subseteq P$ if and only if $K\subset P$ which is expressed 
by $2^m$ rational linear inequalities in $d+m$ variables. Therefore, they define a rational
polytope $\cZ({\bf v})$ in $\R^{d+m}$ (the boundedness of $\cZ({\bf v})$ follows from that of $P$).

The function $|\cdot|:\cZ({\bf v})\to\R$, $Z\mapsto |Z|$ is 
determined by the sum of the last $m$ coordinates in $\R^{d+m}$, hence, is an integer linear function on $\cZ({\bf v})$.
\end{pf}

Notice that reordering of the $v_i$ does not change the zonotope $Z$, so the polytope $\cZ({\bf v})$, as well as the
function $|\cdot|:\cZ({\bf v})\to\R$, is invariant under permutations of the last $m$ coordinates.

Now we are ready for our main result.

\begin{Th}\label{T:period} 
Let $P\subset\R^d$ be a lattice polytope. Then
$$\la(P)=\frac{L(kP)}{k},$$
for some $k\in\N$.
\end{Th}

\begin{pf} 
Consider the polytope $tP$ for some $t\in\N$. It follows from \rp{bound}, that $tP$ contains a smallest
maximal decomposition $Z$ with $m\leq M\vcentcolon=2^d-1$ distinct summands
\begin{equation}\label{e:zonotope2}
Z=a+\al_1[0,v_1]+\dots+\al_m[0,v_m],
\end{equation}
where $a\in\Z^d$,  $v_i\in\Z^d$ are primitive, and $\al_i$ are positive integers whose
sum equals the Minkowski length $L(tP)$. Therefore, $P$ contains a rational zonotope
$$Z/t=a/t+(\al_1/t)[0,v_1]+\dots+(\al_{m}/t)[0,v_m]$$
with $|Z/t|=L(tP)/t$. 
Conversely, every rational zonotope $Z$ in $P$  has the form
\begin{equation}\label{e:zonotope3}
Z=a+\al_1[0,v_1]+\dots+\al_m[0,v_m],
\end{equation}
for some $a\in\Q^d$,  primitive $v_i\in\Z^d$, and non-negative rationals $\al_i$.
Then there exists $t\in \N$ such that $tZ$ is a lattice zonotope in $tP$, and so $|Z|\leq L(tP)/t\leq \la(P)$. 
Therefore, $\la(P)$ is the supremum of the function $Z\mapsto |Z|$ on the set of all rational zonotopes
$Z$ contained in $P$.

We will show below that there exist $\delta>0$, independent of $t$, and a finite set of primitive
vectors $V_\delta\subset\Z^d$ satisfying the following property: If $Z$ is a smallest maximal decomposition in $tP$ for some $t\in\N$ then
$$\la(P)-\delta<|Z/t|\leq\la(P)$$
implies that $v_1,\dots,v_m$ lie in $V_\delta$. By \rl{allzonotopes}, $\la(P)$ equals the
maximum of the linear function $Z\mapsto |Z|$ on the union of rational polytopes $\cZ({\bf v})$ over 
all collections ${\bf v}=(v_1,\dots, v_m)\in (V_\delta)^m$, hence, $\la(P)=|Z'|$ for some rational
zonotope $Z'\subset P$. Choose $k\in\N$
such that $kZ'$ is a lattice zonotope in $kP$. Then 
$$\la(P)=\frac{L(kP)}{k},$$
as required.

It remains to prove the existence of $\delta>0$ and $V_\delta$ satisfying the above property.
Denote $\la=\la(P)$, and $\la_n=\la_n(P)$, the $n$-th rational Minkowski length of $P$. 
Let $e\geq 1$ be the smallest integer such that $\la_n=\la(P)$ for all $n\geq e$. Then to find $\la$ it is enough to 
consider only smallest maximal decompositions in $tP$ of dimension at most $e$. 
Let $Z$ be such a decomposition, as in \re{zonotope2}. 

The case $e=1$ is easy --- we set $\delta=\la/2$ and $V_\delta=U_\delta(P)$, as in \rl{uepsilon}, which is a
finite set. In this case $|Z/t|=\al_1/t\leq s_{v_1}(P)<\delta=\la/2$, unless $v_1\in V_\delta$.

If $e>1$ we have
$$\la_1\leq\dots\leq \la_{e-1}<\la_e=\dots=\la_d=\la.$$ 
Set $\delta=(\la-\la_{e-1})/2$ and choose 
$$0<\ve<\min\Big\{\frac{\la-\delta}{M},\frac{\delta}{M-e}\Big\}.$$

If no $v_i$ lies in $U_\ve(P)$ then $\al_i/t\leq s_{v_i}(P)<\ve$, and so $|Z/t|\leq m\ve<\la-\delta$.
Thus, we may assume that $v_i\in U_\ve(P)$ for $1\leq i\leq k$ and $v_i\not\in U_\ve(P)$ for $k<i\leq m$.

First, suppose that $\{v_1,\dots, v_k\}$ spans an $e$-dimensional subspace. By Lemmas \ref{L:volume}
and \ref{L:finite}, there are only finitely many choices for each $v_i$ for $k<i\leq m$. Thus we define
$V_\delta$ to be the union of $U_\ve(P)$ and the finite sets $V(B)$ for every subset  
$B=\{u_1,\dots, u_e\}\subset U_\ve(P)$ which spans an $e$-dimensional subspace.

Next, suppose the dimension of the span of $\{v_1,\dots, v_k\}$ is less than $e$. We may assume that
$v_l,\dots, v_m$ lie outside of this span for some $k<l\leq m$. Then we have
$$|Z/t|\leq \la_{e-1}+(\al_{l}/t)+\dots+(\al_m/t)<\la_{e-1}+(m-l)\ve.$$
By the choice of $\ve$, and since $l>e$, the latter is smaller than $\la-\delta$.
\end{pf}

\begin{Rem}\label{R:n-rational} 
The same arguments as above show that for any $1\leq n\leq d$, 
$$\la_n(P)=\frac{L(k_nP)}{k_n},$$
for some $k_n\in\N$. In particular, all $\la_n(P)$ are rational numbers.
\end{Rem}

\subsection{Quasi-linearity of the Minkowski length}
The result of \rt{period} allows us to make the following definition.
\begin{Def}\label{D:period} The smallest $k\in\N$ satisfying $\la(P)=L(kP)/k$ is
called the {\it period} of $P$.
\end{Def}

In \rt{quasi-lin} we prove that the Minkowski length is eventually quasi-linear, 
but first we are going to show that the rational Minkowski length is linear (\rp{lin}).
 We will need the following lemma.

\begin{Lem}\label{L:1}
Let $k$ be the period of $P$. 
Then $L(tP)=t\la(P)$ whenever $k\mid t$, $t\in\N$.
\end{Lem}
\begin{pf} Since $tP$ contains the Minkowski sum of $t/k$ copies of $kP$, we have
 $$L(tP)\geq (t/k)L(kP)= t\lambda(P).$$ 
 On the other hand, $t\la(P)\geq L(tP)$, by the definition of $\la(P)$.
\end{pf}

\begin{Prop}\label{P:lin}
For any $t\in\N$ we have $\lambda(tP)=t\lambda(P)$.
\end{Prop}
\begin{pf} We have
$$\lambda(tP)=\max_{s\in\N}\frac{L(stP)}{s}=t\max_{s\in\N}\frac{L(stP)}{st}\leq t\lambda(P).
$$
On the other hand, by \rt{period} and \rl{1},
$$t\lambda(P)=\frac{kt\la(P)}{k}=\frac{L(ktP)}{k}\leq \lambda(tP),
$$ 
where $k$ is the period of $P$.
\end{pf}

\begin{Th} \label{T:quasi-lin}
Let $P$ be a lattice polytope in $\R^d$ with period $k$. Then 
the function $L(tP)$ is eventually quasi-linear with period at most $k$. 
More explicitly, there exist integers $c_r$, for $0\leq r< k$, such that $L(rP)\leq c_r\leq r\la(P)$ and 
$$L(tP)=k \la(P)\left\lfloor\frac{t}{k}\right\rfloor+c_r$$
whenever $t\equiv r\!\!\mod k$ and $t$ is large enough.
\end{Th}

\begin{pf} Fix $0\leq r< k$ and let $t\equiv r\!\!\mod k$. Denote $c_r(t)=L(tP)-k \la(P)\left\lfloor\frac{t}{k}\right\rfloor$.
Note that $c_r(t)$ is an integer as $k\la(P)=L(kP)$. We will show that $c_r(t)$ is constant for $t$ large enough. 
Indeed, $c_r(t)$  are bounded from above:
$$c_r(t)=L(tP)-k \la(P)\left\lfloor\frac{t}{k}\right\rfloor\leq t \la(P)-k \la(P)\left\lfloor\frac{t}{k}\right\rfloor=r\la(P).$$
Also they increase:
$$c_r(t+k)=L(tP+kP)-k \la(P)\left\lfloor\frac{t+k}{k}\right\rfloor\geq L(tP)+L(kP)-k \la(P)\left\lfloor\frac{t+k}{k}\right\rfloor
=c_r(t),$$
where we used $L(kP)=k\la(P)$ in the last equality. Therefore, the integers  $c_r(t)$ for $t\equiv r\!\!\mod k$ 
eventually stabilize to a constant $c_r$.

We have already seen that $c_r\leq r\la(P)$. For the other inequality, let $t=qk+r$. Then, using \rl{1}, we obtain
$$c_r(t)=L(tP)-qk\la(P)\geq L(qkP)+L(rP)-qk\la(P)=L(rP).$$

\end{pf}

\begin{Rem}
The above proof works just as well if we replace $L(P)$ with $L_n(P)$ for any $1\leq n\leq d$ and apply
\rr{n-rational}. Therefore, each $n$-th Minkowski length of $P$ is eventually quasi-linear with period at
most $k_n$. Since $L_1(P)=\lfloor\la_1(P)\rfloor$, the function $L_1(tP)$ is, in fact, quasi-linear.
\end{Rem}

\section{Dimension two}\label{S:dim2}

In this section we deal with lattice polytopes in dimension two. We prove an upper bound 
on the rational length of $P$ in terms of other well-known invariants of $P$ --- the 
Euclidean area and the lattice width of $P$. As an application we give a formula
for $\la(tP)$ and $L(tP)$ for any triangle $P$ in $\R^2$. 

Let $P\subset\R^d$ be a lattice polytope and $v\in\Z^2$ a primitive vector. Recall that the {\it lattice width
of $P$ in the direction $v$} is the integer
$\max_{x\in P}\langle x,v\rangle-\min_{x\in P}\langle x,v\rangle$. Here $\langle x,v\rangle$
is the standard inner product in $\R^d$. The smallest lattice width over all primitive $v\in\Z^d$
is called the {\it lattice width of $P$} and is denoted $w(P)$.

\begin{Prop}\label{P:polygon} 
Let $P\subset\R^2$ be a lattice polygon. Then $\la(P)\leq \frac{2\Vol_2(P)}{w(P)}$
where $\Vol_2(P)$ is the Euclidean area and $w(P)$ the lattice width of $P$.
\end{Prop}

\begin{pf} By the proof of \rt{period}, $\la(P)= |Z|$ for some rational zonotope $Z\subseteq P$ 
with at most 3 distinct summands, i.e.
$$Z=a+\al_1[0,v_1]+\al_2[0,v_2]+\al_3[0,v_3],$$
where $v_i\in\Z^2$ are distinct primitive vectors, $a\in\Q^2$, and $\al_i\in\Q$. We have
$|Z|=\al_1+\al_2+\al_3$.

Let $w_i$ be the lattice width of $P$ in the direction of a primitive vector $v_i^\perp$, orthogonal to $v_i$. We claim that
\begin{equation}\label{e:polygon}
\al_1w_1+\al_2w_2+\al_3w_3\leq 2\Vol_2(P).
\end{equation}
Indeed, let $A_i$ (resp. $B_i$) be a vertex of $P$ where the inner product
 $\langle x,v_i^\perp\rangle $ attains its minimum (resp. maximum). Similarly, 
 let $E_i$ (resp. $I_i$) be the side of $Z$ where  $\langle x,v_i^\perp\rangle $ attains its minimum 
 (resp. maximum). Connect $A_i$ to $E_i$ and $B_i$ to $I_i$ for $i=1,2,3$ by line segments. Also
 triangulate $Z$ (if it is not one-dimensional) by drawing the diagonals through the center of $Z$.
 We obtain a (not necessarily convex) triangulated polygon $S$  inside $P$, see \rf{polygon}.
 
\begin{figure}[h]
\includegraphics[scale=1]{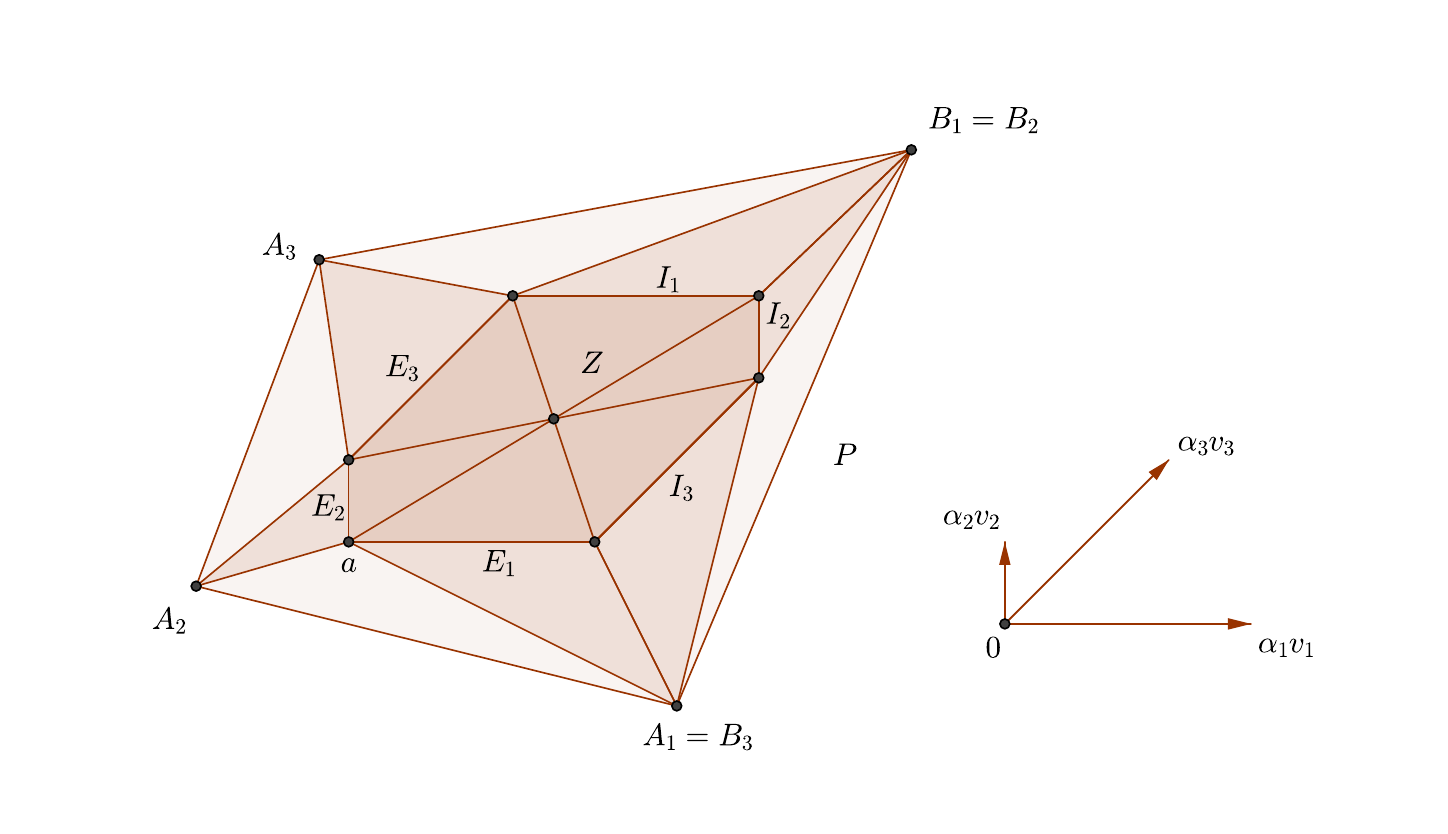}
\caption{A triangulated polygon inside $P$.}
\label{F:polygon}
\end{figure}

 Note that the sum of the areas of the four triangles with bases $E_i$ and $I_i$
 equals $\frac{1}{2}\al_iw_i$.
 Therefore, the left hand side of \re{polygon} represents twice the area of $S$, and \re{polygon} follows.
 
It remains to note that $\al_1w_1+\al_2w_2+\al_3w_3\geq (\al_1+\al_2+\al_3)w(P)=\la(P)w(P)$.
\end{pf} 

Below we apply this bound to give an explicit formula for the (rational) Minkowski length 
of any triangle. Recall the {\it lattice diameter} $\ell(P)$ and the {\it rational diameter} $s(P)$
defined in \rs{lambda}. 

\begin{Cor}\label{C:triangle}
Let $T\subset\R^2$ be a lattice triangle. Let $s(T)$ be its rational diameter and
$\ell(T)$ its lattice diameter. Then 
$$\la(T)=s(T)\quad\text{ and }\quad L(T)=\lfloor s(T)\rfloor =\ell(T).$$
Consequently, $\la(tT)=s(T)t$ and $L(tT)=\lfloor s(T)t\rfloor$. 
\end{Cor}

\begin{pf}
Let $v\in\Z^2$ be a primitive vector such that  the lattice width of $T$ in the 
direction orthogonal to $v$ equals $w(T)$. Then $s_v(T)w(T)=2\Vol_2(T)$, where $s_v(T)$ as in \rd{diameter}.
It follows that $s_v(T)$ is, in fact, $s(T)$. Applying \rp{polygon}, we get  $\la(T)\leq s(T)$.
Conversely, $T$ contains the segment $E$ parallel to $[0,s(T)v]$. Therefore, by the proof of \rt{period},
$\la(T)\geq |E|=s(T)$.

As for $L(T)$, it is clear that $L(T)\leq\la(T)=s(T)$, by the definition of $\la(P)$. 
Thus, $L(T)\leq\lfloor s(T)\rfloor$. On the other hand, $T$ contains a translation of the lattice segment
$[0,\lfloor s(T)\rfloor v]$, hence $L(T)\geq \lfloor s(T)\rfloor$.

Finally, by above, $\la(tT)=s(tT)=ts(T)$ and $L(tT)= \lfloor s(tT)\rfloor= \lfloor ts(T)\rfloor$.
\end{pf}

\begin{Rem}\label{R:tight} The above proof shows that our bound in \rp{polygon} is tight, as 
$\lambda(T)=s(T)=2\Vol_2(T)/w(T)$ for any lattice triangle $T$.
\end{Rem}

\section{Examples and open problems}

In this section we illustrate our results with several examples and raise some questions.

Our first example shows that $L(tP)$ can have an arbitrarily large period $k$.

\begin{Ex} Let $T_k\subset\R^2$ denote the triangle with vertices $(0,0)$, $(k,1)$, and $(1,k)$, for $k\geq 2$.
It is not hard to see that $\ell(T_k)=k-1$ and $s(T_k)=k-1/k$. By \rc{triangle}, 
$$L(tT_k)=\left\lfloor\Big(k-\frac{1}{k}\Big)t\right\rfloor,$$
which is a quasi-linear function with period $k$.
\end{Ex}

\begin{Ex} Let $P$ be a square with vertices $(2,0)$, $(3,2)$, $(1,3)$, and $(0,1)$. One readily sees that 
$\Vol_2(P)=5$, $w(P)=3$, and $L(P)=3$. Therefore, by \rp{polygon}, $\la(P)\leq 10/3$. 

Note that $3P$ contains a
zonotope $Z$ with $|Z|=10$ (in fact, $Z$ is a square with vertices $(2,2)$, $(7,2)$, $(2,7)$, and $(7,7)$).
Therefore, $10\leq L(3P)\leq 3\la(P)\leq 10$, which implies that $\la(P)=10/3$ and $L(tP)$ has period $k=3$.
By \rl{1}, if $t=3q$ then $L(tP)=10q$. Now if  $t=3q+1$ we have
$$10q+3=L(3qP)+L(P)\leq L(tP)\leq 10t/3=10q+10/3,$$
and, hence, $L(tP)=10q+3$. Similarly, $L(tP)=10q+6$ when $t=3q+2$. Therefore,
$$L(tP)=10\left\lfloor\frac{t}{3}\right\rfloor+3r,\ \ \text{ if }\  t\equiv r\!\!\!\!\mod 3.$$
\end{Ex}

Looking at the above examples one may suspect that $L(P)=\lfloor\la(P)\rfloor$
for any polytope $P$. This would imply that $L(tP)$ is not just eventually quasi-linear,
but quasi-linear:
$$L(tP)=k \la(P)\left\lfloor\frac{t}{k}\right\rfloor+L(rP),$$
for any $t\equiv r\!\!\mod k$ (see \rt{quasi-lin}). 

However, $L(P)=\lfloor\la(P)\rfloor$ does not hold even for the case of lattice polygons, as
demonstrated by the following example.

\begin{Ex}
Let $P=2Q$ where $Q$ is the square with vertices $(1,0)$, $(5,1)$, $(4,5)$, and $(0,4)$.
Then by \rp{polygon},
$$\la(P)=2\la(Q)\leq\frac{68}{5}.$$
Also  $L(5P)$ contains $Z$ with $|Z|=68$ (namely, $Z$ is the square with vertices $(8,8)$, $(8,42)$, $(42,42)$, and $(42,8)$), hence, $\la(P)=68/5$.  By observation, $L(P)=12$, which illustrates that $\la(P)-L(P)$ can be as large as $8/5$.
\end{Ex}

\begin{Prob} Find the supremum of $\la(P)-L(P)$ over all lattice polytopes $P\subset\R^d$.
\end{Prob}

It is not hard to see that $\la(P)-L(P)<4$ for any lattice polygon $P\subset\R^2$, but we are confident that this
bound could be improved.

In all 2-dimensional examples we computed, the function $L(tP)$ was always quasi-linear.
Although we do not expect this to be the case in general, we have not been able to produce a counterexample.

\begin{Prob} Prove that $L(tP)$ is quasi-linear or give an example of a lattice polytope $P$ for which $L(tP)$ is not quasi-linear.
\end{Prob}

Finally, we have seen that $L(P)=L_1(P)$ when $P$ is a positive integer dilate of a unimodular simplex as well as when $P$ 
is any simplex in dimension two. This prompts the following problem.

\begin{Prob} Prove or disprove  that for any simplex in $\R^d$ its Minkowski length coincides with its lattice diameter.
\end{Prob}

\end{document}